\documentclass{article}
\usepackage[T1]{fontenc}
\usepackage[utf8]{inputenc}
\usepackage{authblk}
\usepackage{amsmath,amssymb,amsfonts,amsthm,amsbsy}
\usepackage[dvips]{graphicx}
\usepackage{subfigure}
\usepackage{epsfig}

\newtheorem{THM}{Theorem}[section]

\newtheorem{LEM}[THM]{Lemma}
\newtheorem{PRP}[THM]{Proposition}
\newtheorem{DEF}[THM]{Definition}

\newcommand{\nnb}{\nonumber}

\newcommand{\E} {{\mathbb E}}
\renewcommand{\P} {{\mathbb P}}
\newcommand{\Q} {{\mathbb Q}}

\newcommand{\Zbold} {{\mathbb Z}}
\newcommand{\Nbold} {{\mathbb N}}

\newcommand{\Gcal} {{\cal G}}

\newcommand{\Ocal} {{\cal O}}

\newcommand{\ind}{\mathbf{1}}

\newcommand \nsde {\alpha}

\title{Conditioning the logistic branching process on non-extinction}

\author[1]{Alison Etheridge\thanks{Supported in part by EPSRC Grant EP/101361X/1}}
\author[1]{Shidong Wang\thanks{Supported in full by EPSRC Grant EP/101361X/1}}
\author[2]{Feng Yu\thanks{Supported in part by EPSRC Grant EP/I028498/1}}
\affil[1]{Department of Statistics, 1 South Parks Road, University of Oxford, Oxford, OX1 3TG, UK}
\affil[2]{School of Mathematics, University Walk, University of Bristol, Bristol, BS8 3AN, UK}

\begin{document}
\maketitle
\pagestyle{headings}

\begin{abstract}
We consider a birth and death process in which death is due to both `natural death'
and to competition between individuals, modelled as a quadratic function of population size.
The resulting `logistic branching process' has been proposed as a model for numbers of individuals
in populations competing for some resource, or for numbers of species.  However, because of the
quadratic death rate, even if the intrinsic growth rate is positive, the 
population will, with probability one, die out in finite time. There is considerable 
interest in understanding the process conditioned on non-extinction.

In this paper, we exploit a connection with the ancestral selection graph of population
genetics to find expressions for the transition rates in the logistic branching process conditioned on survival
until some fixed time $T$, in terms of the distribution of a certain one-dimensional 
diffusion process at time $T$.  We also find the probability generating function of
the Yaglom distribution of the process and rather explicit expressions for the 
transition rates for the so-called Q-process, that is the logistic branching process conditioned
to stay alive into the indefinite future. For this process, one can write down the joint generator of the (time-reversed)
total population size and what in population genetics would be called the `genealogy' and in phylogenetics would
be called the `reconstructed tree' of a sample from the population. 

We explore some ramifications of these calculations numerically.
\end{abstract}

\section{Introduction}

\subsection{Background}

Models of population growth are of central importance in mathematical ecology.  Their origin can 
be traced back at least as far as the famous essay of Thomas Malthus in 1798. \nocite{malthus:1798}
Probably the most popular stochastic models are the classical Galton-Watson branching processes, or
their continuous time counterparts, the Bellman-Harris processes.  The key assumption that they make, that
all individuals
in the population reproduce independently of one another, is extremely convenient mathematically.
However, the elegance and tractability of the branching process model is offset by the difficulty that
it predicts that the population will, with probability one, either die out in finite time
or grow without bound: one would prefer
a model that predicted more stable dynamics.
 
In population genetics, one typically circumvents this problem by conditioning
the total population size to be identically constant.
But, as argued for example by Lambert~(2010),\nocite{lambert:2010} 
a much more satisfactory solution would
be to take a `bottom-up' approach and begin with an individual based model.  

An alternative to conditioning on constant population size is 
to observe that since we are able to sample from the population, we are necessarily observing a 
realisation of the population process conditioned on non-extinction. The difficulty then is that one needs
additional information, such as the age of the population, to know how it has evolved to its current state.  
Moreover, as the age of the population grows, unless the underlying branching process is subcritical, the size
of the population conditioned on survival
grows without bound and so this does not present a good model for the relatively stable 
population sizes that we often observe (at least over the timescales of interest to us) in the wild.
Again as argued by
Lambert~(2010),\nocite{lambert:2010}
even if it is doomed to ultimate extinction,
the size of an isolated population can fluctuate for a very long time relative to our chosen timescale, and
then it makes sense to consider the dynamics of a population conditioned to survive indefinitely long into the
future.  This is the so-called Q-process, which we define more carefully in Definition~\ref{defn:Q-process}.

One difficulty with the branching process model is that it does not impose any restriction on the
population size.  In reality, we expect that as population size grows, competition for resources will
reduce the reproductive success of individuals and/or their viability.  
Our goal here is to study a model which introduces this effect in the simplest possible way.  That is,
we shall consider the birth and death process described, for example, in Chapter 11, Section 1, Example A of
Ethier \& Kurtz~(1986),
\nocite{ethier/kurtz:1986}
in which births and `natural deaths' both occur at rates proportional to the current population size (as in the
classical birth and death process), but there are also additional deaths, due to competition, that occur at a rate
proportional to the square of the population size.  This quadratic death rate prevents the population from
growing without bound, but even if the `intrinsic growth rate' determined by the births and natural deaths is
positive, the population will, with probability one, die out in finite time.

Let us give a precise definition of the stochastic process that we shall study.
\begin{DEF}
\label{defn:logistic branching}
The logistic branching process, $Z^{b,c,d}=(Z^{b,c,d}_t,t\ge 0)$,
is a population model in which each individual gives birth at rate $b$,
dies naturally at rate $d$, and dies due to competition at rate  $c(i-1)$ where $i$ is the
current population size.
In particular, it is a pure jump process
taking values in $\Zbold^+$, with jump rates $q_{ij}$
\[ q_{ij} = \left\{ \begin{array}{ll}
  di+ci(i-1), & \mbox{if $i\ge 1$ and $j=i-1$} \\
  bi, & \mbox{if $i\ge 1$ and $j=i+1$} \\
  0, & \mbox{otherwise}
\end{array} \right. .
\]
\label{def0}
\end{DEF}
For definiteness, in our mathematical arguments, 
we shall use the language of population models.  Thus the logistic branching process models the size of
a stochastically fluctuating population.  
However, in the applications we have in mind,
the `individuals' in the population may be species.
If $c=0$, then we recover the classical birth and death process.

Our aim in this paper is to consider the logistic branching process conditioned on non-extinction.
Because the quadratic death rate destroys the independent reproduction which made the classical branching
process so tractable, we have to work much harder to establish the transition rates of our conditioned 
process.  We shall establish these both for the process conditioned to stay alive until some fixed
time $T$ and for the Q-process, with those for the Q-process being rather explicit.  

\subsection{Reconstructed trees}
\label{reconstructed}

There are a number of interesting objects to study within this model.  
For example, in population genetics, one is interested in the genealogical 
trees relating individuals in a sample of individuals from the population. 
This will be a random timechange of the classical Kingman coalescent, but 
the timechange is determined by the distribution of the path of the population
size as we trace backwards in time. Since our Q-process is reversible, we can 
explicitly write down the joint generator of the coalescent and the total population
size in that setting (Theorem~\ref{thm:coalescent}).  
In the context of phylogenetics, one is 
also interested in this 
object (where the sample may be the whole population).  Individuals then
correspond to species, and the genealogy of the population is usually called the 
`reconstructed tree', Nee et al.~(1994).\nocite{nee/may/harvey:1994} 
It corresponds to the phylogenetic tree in which all lineages that have
terminated by the present time have been removed. 
There is a long tradition in paleontology of using simple mathematical models as tools
for understanding patterns of diversity through time.  Nee~(2004) provides a 
\nocite{nee:2004}
brief survey, emphasizing the predominance of analyses based on birth and death processes
and on the Moran process. In this second model, which is familiar from population 
genetics, the total number of lineages
remains constant through time: the extinction of a lineage is matched by the birth
of a new lineage.  Our model is in some sense intermediate between these two classes. 

The reconstructed tree constructed from a standard birth-death process model will exhibit what has been dubbed the 
`pull of the present', Nee et al.~(1992). \nocite{nee/mooers/harvey:1992}
It arises from the fact that lineages arising in the recent past are more likely to be 
represented in the phylogeny at the present time
than lineages arising in the more distant past, simply because they
have had less time to go extinct.  
By looking at the way lineages have accumulated through time, one is able to estimate both the birth and
death rates from the reconstructed tree.  However, real phylogenies rarely exhibit the rate of accumulation 
of lineages through time predicted by a birth-death model.  In particular, one sees a `slowdown' towards
the present, see e.g.~Etienne \& Rosindell~(2012) and references therein. 
\nocite{etienne/rosindell:2012}
They propose a `protracted speciation' model, in which speciation is divided
into two phases: `incipient' and `good'. An incipient
species does not branch and produce new species. An alternative explanation
of the apparent slowdown in diversification can be found in Purvis {\em et al.}~(2009), who
\nocite{purvis/etal:2009}
argue that it is simply
due to age-dependency in whether or not nodes are deemed to be speciation events.
A simple way to model this is to attach an exponential clock with
rate $\lambda$ to each species:
only species older than the corresponding exponential random variable are `detectable', and only detectable species
can be sampled, whereas both detectable and undetectable species can branch
and produce new species. 

A commonly used metric to determine the shape of a phylogeny is
the $\gamma$ statistic. It is based on $g_2, g_3,\ldots,g_n$,
the internode distances of a reconstructed phylogeny with $n$ taxa, for example
$g_2$ is the length of the time period during which there are 
exactly 2 lineages in the reconstructed tree.
The $\gamma$ statistic is defined as (Pybus \& Harvey~2000)
\nocite{pybus/harvey:2000}
\[ \gamma = \frac{(\frac 1 {n-2} \sum_{i=2}^{n-1} \sum_{k=2}^i k g_k)
  - \frac 1 2 \sum_{j=2}^{n} j g_j}{
  \sum_{j=2}^{n} j g_j \sqrt{\frac 1 {12(n-2)}}}
\]
Under a pure birth process, the expected value of $\gamma$ is 0.
If $\gamma<0$, then the reconstructed tree's internal nodes are closer to the root than
the tip, and {\em vice versa} for $\gamma>0$. 

Phillimore \& Price~(2008) measured
the $\gamma$ statistic for 45 clades of birds and obtain $\gamma$ values
ranging from -3.26 to 1.85.
\nocite{phillimore/price:2008}
They argue 
that their data is strongly suggestive of density-dependent 
speciation in birds.  They predict that 
the rate of speciation slows down as ecological opportunities and
geographical space place limits on the opportunities for new species to develop.
This differs from our model, in which it is the death rate rather than the birth
rate of species that is density dependent. However, the difference is shortlived, in
the sense that species born under our model which would not have appeared if the
birth rate were density dependent, are rapidly removed by the density-dependent death
rate.  One might expect that if we only consider `detectable' species, in the
sense of Purvis {\em et al.}~(2009), then we should recover something akin to
a model with density-dependent birth rate.  

\begin{figure}
\centering
\includegraphics[height=0.4\textheight,width=0.8\textwidth]{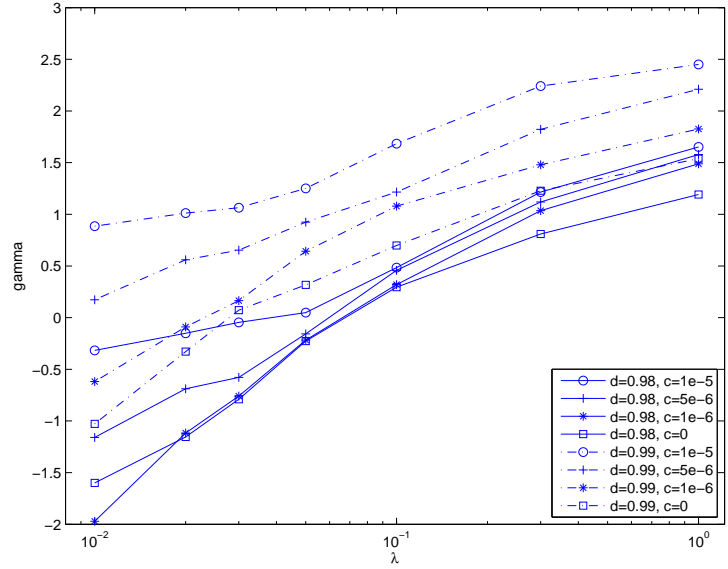}
\caption{Gamma value vs $\lambda$: $b=1$.}
\label{fig:gamma1}
\end{figure}

A plot of the $\gamma$ statistic versus $\lambda$ obtained
from simulating our logistic branching model is shown in
Figure~\ref{fig:gamma1}. It suggests that our logistic branching model, with
birth rate scaled to be one, 
produces realistic values of $\gamma$ if $\lambda$ is taken to be $<0.1$.
Thus the logistic branching model together with the effect of delayed
detection is also consistent with the apparent slowdown of accumulation
of lineages.

\subsection{Conditioning on non-extinction}

Because of the quadratic death rate, it is clear, as observed for example by Lambert~(2005),
\nocite{lambert:2005}
that a population evolving according to the model in Definition~\ref{def0} will die out in finite time.
We shall circumvent this problem by conditioning on non-extinction and looking at
stationary behaviour of the conditioned process.
When 0 (or any other state) is an absorbing state of the process,
rendering a non-trivial stationary distribution impossible,
there are usually two `quasi-stationary' objects one can study. The first one corresponds to 
conditioning the process on non-extinction at the present time.
\begin{DEF}[Quasi-stationary distributions]
\label{defn:quasi-stationary distributions}
We say that the probability measure $\nu$ is a quasi-stationary distribution of the process $\{Z_t\}_{t\geq 0}$ if 
$$\nu(\cdot)=\P_\nu(Z_t\in \cdot \ | Z_t>0).$$
\end{DEF}
Of particular importance is the Yaglom limit:
\begin{DEF}
\label{defn:yaglom limit}
The Yaglom limit, if it exists, is the probability measure $\nu$ obtained as 
$$\lim_{t\to\infty} \P(Z_t/g(t)\in \cdot \ | Z_t>0),$$
where $g(t)$ is a deterministic function.
\end{DEF}
In our setting, because the quadratic competition controls the size of the
conditioned population,
no normalising function is required and we may set $g(t)\equiv 1$ (see below). 
However, although we are able to find expressions both for the probability 
generating function of the Yaglom limit and the
transition rates of $Z_u|Z_t\neq 0$, for $0<u<t$,
neither is in a particularly convenient form.

The second notion of conditioning on non-extinction, which turns out to be somewhat
simpler in our
setting, corresponds to conditioning the process to 
survive into the indefinite future.  This yields the so-called Q-process.
\begin{DEF}[Q-process]
\label{defn:Q-process}
Let $\tau_0$ be the hitting time of the absorbing state $0$.
The Q-process corresponding to $\{Z_t\}_{t\geq 0}$ is determined as follows. 
Write $\{{\cal F}_t\}_{t\geq 0}$ for the natural filtration, then for any $t>0$ and
any ${\cal F}_t$-measurable set $B$,
\[ \Q(Z\in B) = \lim_{u\to\infty} \P(Z\in B \ | \ \tau_0>t+u). \]
\end{DEF}
Our main analytic result is the following.  
\begin{THM}
Let us write $Z^*$ for the Q-process derived from the logistic branching process $Z^{s,\nsde,\mu}$
of Definition~\ref{defn:logistic branching} and
$\{q^*_{i,j}\}_{i,j\in \Nbold}$ for the corresponding transition rates. Then $q^*_{i,j}=0$ for $|i-j|>1$ and
\[ q^*_{k,k+1} = q_{k,k+1} r^*_{k+1,k}, \  \forall k\geq 1;\qquad 
  q^*_{k,k-1} = q_{k,k-1} r^*_{k-1,k},\ \forall k\geq 2, \]
where
\begin{eqnarray}
  r^*_{i,j} = \frac{\int_0^1 (1-(1-\zeta)^i) \pi^*(\zeta) \ d\zeta}{
    \int_0^1 (1-(1-\zeta)^j) \pi^*(\zeta) \ d\zeta},
\label{def:rstar}
\end{eqnarray}
and $\pi^*$ is defined in Proposition~\ref{prp:pi_star}.
Moreover, there exists $C>0$, independent of $k$, such that
$$\frac{1}{C}<r^*_{k-1,k}\leq r^*_{k+1,k}<C.$$
\label{cor:q_star}
\end{THM}
In the statement of the result we have set $b=s, c=\nsde, d=\mu$.  The reason for these slightly strange looking 
choices will become clear from the proof.

A trivial coupling argument guarantees that $Z_t|Z_{t+u}>0$ stochastically dominates 
$Z_t|Z_t>0$ and so the existence of a non-trivial stationary distribution for the Q-process is
enough to guarantee that we can take $g(t)\equiv 1$ in defining the Yaglom limit for
the logistic branching process.

\subsection{Previous mathematical work}

There is a very substantial mathematical literature devoted to population processes with
density dependent regulation.  Most are considerably more complex than that proposed here.
For example, (spatial and non-spatial) branching processes with mean offspring number
chosen to depend on (local) population density have been studied by many authors including
Bolker \& Pacala~(1997), Campbell~(2003), Law~{\em et al.}~(2003), Etheridge~(2004).
\nocite{bolker/pacala:1997}\nocite{campbell:2003}\nocite{etheridge:2004}
\nocite{law/murrell/dieckmann:2003}
The model considered here 
is the simplest possible model for a density dependent population
process and as a result we are able to obtain more precise results than have been found for the
more complicated regulated branching processes.  
Moreover, one expects that the qualitative behaviour of our model should mirror that of more complex models.

Previous studies of the logistic branching process include
Lambert~(2005) and Lambert~(2008) who 
\nocite{lambert:2005}\nocite{lambert:2008}
considered both the individual based model of
Definition~\ref{def0}, and the continuous state branching process, sometimes called the Feller diffusion
with logistic growth, which arises as a scaling limit. It will be clear that we could take the 
corresponding scaling limit in our results and we indicate the appropriate scaling in \S\ref{scaling}.  
However, our primary purpose here is to consider the discrete model.

Pardoux \& Wakolbinger~(2011) and Le {\em et al}~(2013)
\nocite{pardoux/wakolbinger:2011}\nocite{le/pardoux/wakolbinger:2013}
study a process with the same transition rates as our logistic branching process,
but in contrast to our setting, individuals in the population are not
exchangeable.  Instead, each has a label and the chance of being killed due to 
competition depends upon that label.  As a result the genealogy of the population in
their model is quite different from the one that interests us in our
biological applications. Their motivation is also
quite different from ours.  They rescale and pass to a diffusion limit and in the process
recover an analogue of the Ray-Knight Theorem for the Feller diffusion with logistic growth.  

Classical studies of quasi-stationary distributions for Markov chains with
finitely many states (Darroch \& Seneta, 1965) and infinitely many
\nocite{darroch/seneta:1965}\nocite{seneta/vere-jones:1966}
states (Seneta \& Vere-Jones, 1966) rely on the Perron-Frobenius Theorem
and finding expressions for the left and right eigenvectors corresponding
to the largest eigenvalue of the transition matrix of the chain.
The logistic branching process we consider
here is a continuous-time Markov chain with infinitely many states, but
the presence of the competition term, makes it difficult to characterise
the left and right eigenvectors of the transition matrix and so we adopt a different approach.

In the continuous setting, 
a study of conditioned diffusion models in population
dynamics was carried out in Cattiaux {\em et al}~(2009), who consider 
\nocite{cattiauxetal:2009}
one-dimensional diffusion processes of the form
\begin{equation}
\label{cattiaux sde}
dX = dB - q(X) \ dt. 
\end{equation}
They establish conditions under which there is a unique quasi-stationary
distribution as well as the existence of the Q-process.
We note that, when properly rescaled, the logistic branching process of
Definition~\ref{def0} converges to a diffusion of the form
$$dZ=\sqrt{\nsde Z}\ dB+h(Z)\ dt.$$
It was equations of this form that motivated the study of Cattiaux {\em et al.}~(2009);
if one defines $X=2\sqrt{Z/\nsde}$ then $X$ satisfies~(\ref{cattiaux sde}) with a suitable
choice of $q$.

The $Q$-process corresponding to a critical or subcritical branching process
can be viewed as an immortal `backbone' which constantly throws off subfamilies which
each die out in finite time.  In the continuous state branching process limit, these 
families each evolve as an independent copy of the unconditioned process and we recover
the immortal particle representation of Evans~(1993). \nocite{evans:1993} 
In the logistic branching model, the Q-process can also be decomposed in this way.
The existence of a unique `immortal particle' follows readily as the state $1$ is
recurrent for the process, even when the `intrinsic growth rate', $b-d$ is positive.  In particular,  
the time that we must trace back before the 
present before we reach the most recent common ancestor of the 
current population, which is certainly dominated by the time since the 
population size was last $1$, is necessarily finite. In the case of the Q-process for a subcritical
continuous state branching process, Chen \& Delmas~(2012) examined the size of 
\nocite{chen/delmas:2012}
the population at the time of that most recent common ancestor.  They found
a mild bottleneck effect: the size of the population just before the MRCA is 
stochastically smaller than that of the current population.

\begin{figure}
\centering
\includegraphics[height=0.4\textheight,width=0.8\textwidth]{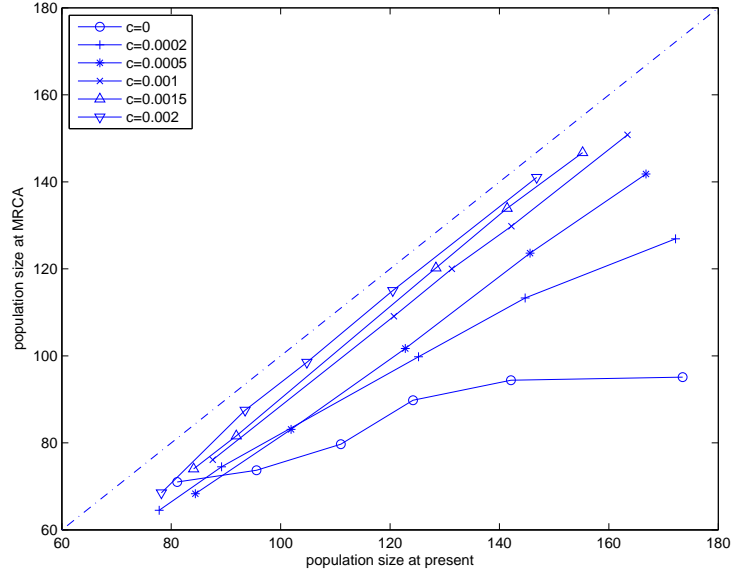}
\caption{Population size at MRCA vs at present for difference values of $c$:
$d=1$, $b\le 1$ varies along each line.}
\label{fig:pop_size_mrca}
\end{figure}

\begin{figure}[!ht]
\centering
\includegraphics[height=0.2\textheight,width=0.8\textwidth]{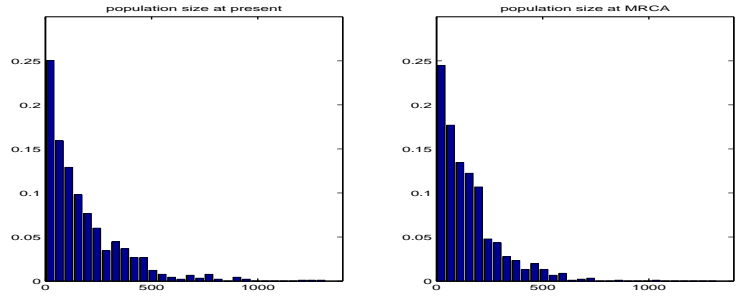}
\caption{Distribution of population size at sampling time $T=500$:
$b=0.995, d=1, c=0$.}
\label{fig:pop_size_comp1}
\end{figure}

\begin{figure}[!ht]
\centering
\includegraphics[height=0.2\textheight,width=0.8\textwidth]{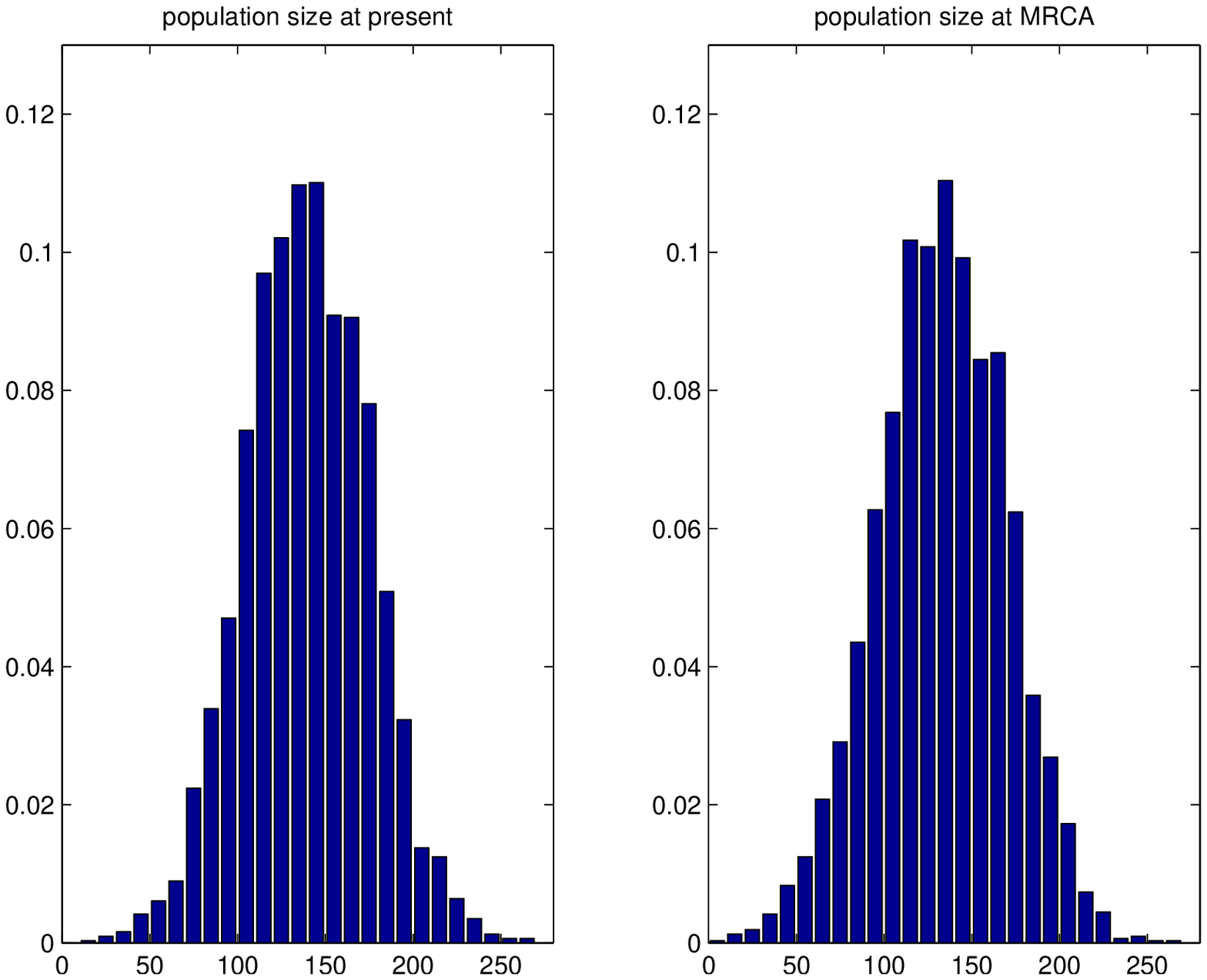}
\caption{Distribution of population size at sampling time $T=500$:
$b=1.15, d=1, c=0.001$.}
\label{fig:pop_size_comp2}
\end{figure}

There is, of course, no reason to expect a similar effect in our model once the intrinsic
growth rate is positive.  Perhaps more surprising is that even in the case where
the intrinsic growth rate is negative, which should more closely resemble the setting
of Chen \& Delmas~(2012), the bottleneck
effect seems to dwindle with the introduction of
competition between individuals. Figure~\ref{fig:pop_size_mrca} shows that
the bottleneck effect becomes less severe as $c$ increases.
Figures~\ref{fig:pop_size_comp1} and~\ref{fig:pop_size_comp2} contrast
the distribution functions of the population sizes for a subcritical branching process
and a logistic branching process. The parameters chosen for
Figures~\ref{fig:pop_size_comp1} and~\ref{fig:pop_size_comp2} are such that
the population sizes at sampling time rough match each other.
In both cases, we plot an approximation of the Yaglom limit, rather than an approximation
of the stationary distribution of the Q-process.
The distribution
of population sizes for the sub-critical branching process is exponential-like,
whereas the distribution for the logistic branching process is Gaussian-like.
We show in \S\ref{sec:weak_comp} that the stationary distribution of
the Q-process of the logistic branching process with weak competition is also
approximately Gaussian.

\subsection{Outline}

The powerful tools available for studying branching processes break down in our setting.
We are also unable to directly exploit the machinery of one-dimensional diffusions.  
However, 
it turns out that we can access that machinery indirectly, through a duality between
our logistic branching process and a Wright-Fisher diffusion process with
selection and one-way mutation. This duality, which is explained in detail in
\S\ref{duality}, arises from interpreting the logistic branching process as the
`ancestral selection graph' of Krone \& Neuhauser~(1997) in a degenerate case (where the
\nocite{krone/neuhauser:1997}
Wright-Fisher diffusion does not have a stationary distribution).
This duality allows us to express the transition rates of the
conditioned logistic branching process in terms of the distribution of the 
Wright-Fisher diffusion conditioned on non-extinction.  In \S\ref{Q-process of diffusion}
we use standard techniques to calculate the stationary distribution of the Q-process of the diffusion
which allows us to write down the transition rates for the conditioned logistic branching 
process in \S\ref{transition rates}.  Our expressions are rather unwieldy for explicit 
calculations and so in \S\ref{perturbation} we find a simple approximation, valid for 
weak competition. 
Finally, in \S\ref{yaglom}, we obtain an expression for the probability generating 
function of the Yaglom limit. 

\section{Construction and properties of the dual SDE}

\subsection{The dual process}
\label{duality}

We will define a stochastic differential equation which is dual, in an 
appropriate sense, to the logistic branching process
defined in Definition~\ref{def0}.  

It is convenient to start from a Moran model
with constant population size $N$.  We write $[N]=\{1,2,\ldots ,N\}$. Each $i\in [N]$ corresponds
to an individual in the population whose genetic type, $X_i$, can take one of two states which
we denote by $a$ and $A$.
We suppose that type $a$ is selectively favoured over type $A$,
with selection coefficient $s$, but there is one-way mutation
from type $a$ to type $A$ at rate $\mu$.
More precisely, the dynamics of the process are driven by independent Poisson processes,
$\Lambda^M_i$, $\Lambda^R_{i,j}$ and $\Lambda^S_{i,j}$
($i, j\in [N]$, $i\ne j$)
with intensities $\mu$, $\nsde/2$ and $s/N$,
respectively. 
The type of each individual, denoted $X_i$, is set to $a$ at $t=0$
and evolves as follows
\begin{enumerate}
\item[1.] Mutation: at each jump of $\Lambda^M_i$, individual $i$ mutates
to type $A$.
\item[2.] Resampling: at each jump of $\Lambda^R_{i,j}$, the type of individual
$i$ is replaced by that of individual $j$.
\item[3.] Selection: at each jump of $\Lambda^S_{i,j}$, if individual $i$
is of type $A$ and individual $j$ is of type $a$ immediately before this jump, then
the type of individual $i$ becomes $a$; otherwise, nothing happens.
\end{enumerate}
Let $P_t$ denote the proportion of type $a$ individuals in the population
at time $t$, then $P_0=1$ and standard techniques
(see, for example, the calculation in the proof of Theorem 10.1.1
in Ethier \& Kurtz 1986) imply that for each fixed $T>0$, as $N\to\infty$,
$\{P_t, t\in [0,T]\}$ converges weakly to $\{p_t, t\in [0,T]\}$,
where $p$ solves the following stochastic differential equation
\begin{eqnarray}
  dp = (-\mu p + sp(1-p)) \ dt + \sqrt{\nsde p(1-p)} \ dW, \ p_0=1.
\label{eq:sde}
\end{eqnarray}
Moreover, we can control the rate of convergence.
\begin{PRP}
\[ \E \left[ \sup_{t\in [0,T]} |P(t)-p(t)| \right] \le C_T N^{-1/2}. \]
\label{prp:sde_conv}
\end{PRP}

Now suppose that we sample $k$ individuals from the population at time $t$.
We construct a system of branching and coalescing lineages,
denoted by $\xi^{s,\nsde, \mu}$ in which the ancestry of our sample is embedded. 
When there is no ambiguity, we drop
the superscripts. 
Time for the process $\xi$ runs backwards from the perspective of the Moran model. 
Thus our sample is taken at time zero for $\xi$ (which is $t$ for the Moran model)
and we trace the ancestry backwards to time $t$ for $\xi$ (which is time $0$ for the
Moran model).  
We use the same labels as for our Moran model, to label the lineages in $\xi$, so
that $\xi_u\subseteq [N]$ for each time $0\leq u\leq t$.
We shall write $\kappa_u$ for the number of lineages in $\xi_u$.   

The evolution of $\xi$ is completely determined by the Poisson processes
$\Lambda^M_i$, $\Lambda^R_{i,j}$ and $\Lambda^S_{i,j}$ that were used
to construct the forwards-in-time process $P$. More specifically,
\begin{enumerate}
\item[1.] Mutation: if $i\in\xi_u$ and there is a jump of
$\Lambda^M_i$ at time $t-u$, then $i$ is removed from $\xi$ at time $u$.
\item[2.] Resampling: if $i,j\in \xi_u$ and there is a jump of
$\Lambda^R_{i,j}$ at time $t-u$, then $i$ is removed from $\xi$ at time $u$.
\item[3.] Selection: if $i\in \xi_u$ and there is a jump of
$\Lambda^S_{i,j}$ at time $t-u$ for any $j\in [N]$,
then $j$ is added to $\xi$ at time $u$ (if it is not already a member).
\end{enumerate}
The construction of the $\xi$ here closely resembles that of the ancestral
selection graph of Krone \& Neuhauser~(1997). In particular, if we assign types
\nocite{krone/neuhauser:1997}
to the lineages alive at time $t$ for $\xi$ ($0$ for the Moran model), then we can work
our way
through the branching and coalescing structure $\xi$ to deduce the types of individuals
in the sample from the population at time $0$ ($t$ for the Moran model). 
The main modification of the usual ancestral selection graph is that,
since mutation is always from $a$ to $A$, it is
not necessary to trace any lineage beyond the first mutation that we encounter - we
have already discovered that it must be of type $A$.
The effects of the resampling and selection mechanisms
on $\xi$ are exactly the same as in the ancestral selection graph: 
a selection event falling on lineage $i$ causes it to branch,
producing another lineage $j$. In order to know which of the two resultant lineages gives its type
to the lineage we are tracing, we need to know the types of both lineages $i$ and $j$ immediately
after (before in Moran time) the selection event.

The number of lineages $\kappa^{s,\nsde, \mu}$ in the dual process
$\xi^{s,\nsde, \mu}$
evolves backwards in time according to the following jump rates:
\begin{enumerate}
\item[1.] Mutation: $\kappa$ is decremented by 1 at rate $\mu \kappa$.
\item[2.] Resampling: $\kappa$ is decremented by 1
  at rate $\nsde \kappa(\kappa-1)$.
\item[3.] Selection: $\kappa$ is incremented by 1
  at rate $s \kappa (1-\frac{\kappa-1}{N-1})$.
\end{enumerate}
Thus as $N\to\infty$, the evolution of the genealogical process $\xi_{\cdot}$
is almost the same as that of the logistic branching process.

\begin{LEM}
For any $T>0$, there exists a constant $c_T$, independent of $N$, such that
the paths of processes $Z^{s,\nsde,\mu}_\cdot$ with initial
condition $Z^{s,\nsde,\mu}_0=\kappa_0$ and $\kappa^{s,\nsde,\mu}_{\cdot}$
with initial condition $\kappa^{b,c,d}_0=\kappa_0$
coincide up to time $T$ with probability at least $1-c_T N^{-1/3}$.
\label{lem:couple1}
\end{LEM}

\proof
We couple the processes $Z^{s,\nsde,\mu}$ and $\kappa^{s,\nsde,\mu}$ so that
$\kappa^{s,\nsde,\mu}\le Z^{s,\nsde,\mu}$ a.s. for all $t\ge 0$.
We fix $T$. Let
\[ A=\left\{ \sup_{t\in [0,T]} Z^{b,c,d}_t > K e^{bT} \right\} \]
and $\tau$ be the first time when the paths of
$Z^{s,\nsde,\mu}$ and $\kappa^{s,\nsde,\mu}$ diverge.
Since
\[ \E\left[\sup_{t\in [0,T]} Z^{b,c,d}_t \right] \le C_1 e^{bT}, \]
we have
\[ \P(A) < C_1/K.\]
We start $Z^{s,\nsde,\mu}$ and $\kappa^{s,\nsde,\mu}$ at the
same initial position $\kappa_0$ and restrict our attention to $A$.
We couple $Z^{s,\nsde,\mu}$ and $\kappa^{s,\nsde,\mu}$ until $\tau$.
Since the rates at which
these two processes decrease are exactly the
same, whereas the rates at which they increase differ by $bi(i-1)/(N-1)$,
the probability of $\{\tau<T\}$ is dominated by the probability that
a $Poisson(b K^2 e^{2bT} T/(N-1))$ takes a nonzero value, hence
\[ \P(\{\tau<T\}\cap A^c) \le 1-e^{-b K^2 e^{2bT} T/(N-1)}
  \le b K^2 e^{2bT} T/(N-1). \]
We conclude that
\[ \P(\tau<T) \le \frac{C_1} K + \frac{b K^2 e^{2bT} T}{N-1}. \]
Taking $K=N^{1/3}$ yields the desired result.
\qed

\begin{figure}
\centering
\subfigure[At least one ancestral lineage has reached Moran time $t=0$, hence there is
at least one type $a$ individual in the present sample.]{
\includegraphics[height=0.4\textheight,width=0.8\textwidth]{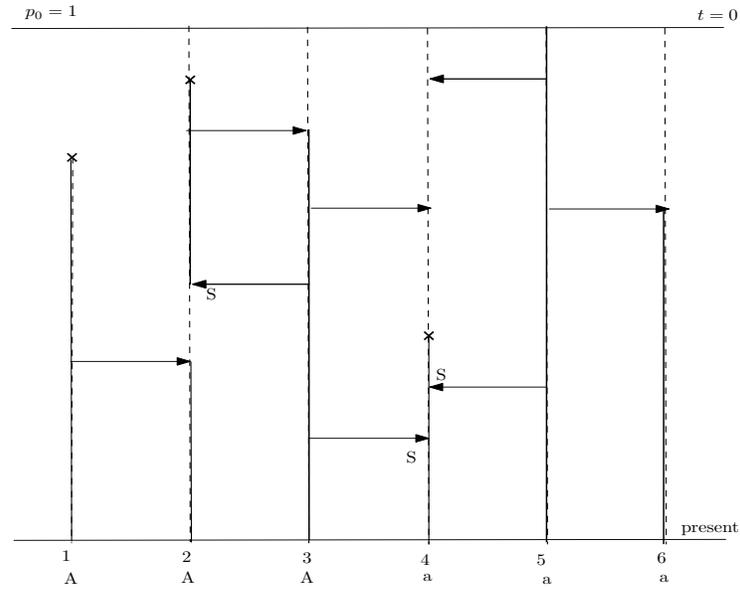}}
\subfigure[No ancestral lineage reaches Moran time $t=0$, hence all individuals
in the present sample are of type $A$.]{
\includegraphics[height=0.4\textheight,width=0.8\textwidth]{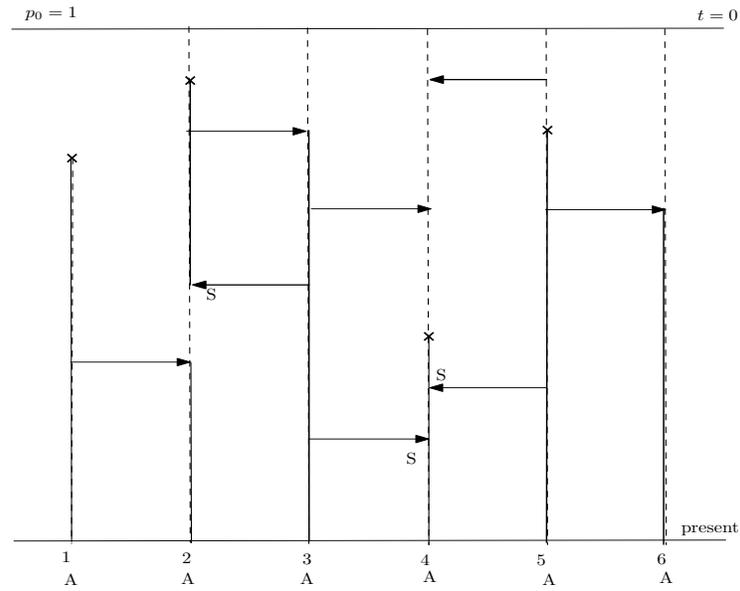}}
\caption{Two realisations of the ancestral selection graph. A cross denotes
a mutation from type $a$ to $A$; an arrow with no letter attached denotes
a resampling event, with the direction of the arrow corresponding to the
direction of replacement; an arrow with S attached denotes a selection event.} 
\label{fig1}
\end{figure}

We take a sample of individuals $\xi_0=\{1,\ldots,\kappa_0\}$
from the population at (Moran) time $t$ and trace back their
lineages through $\xi_{\cdot}$. We refer to Figure~\ref{fig1}
for an illustration of the two possible realisations of the ancestral process.
If $\xi_t\neq\emptyset$, then at least one
lineage has survived until time $t$. By assumption, every individual at that 
time is
of type $a$, and so its descendants will `win' every selection event that they
encounter since type $a$ is
selectively favoured. Thus tracing back to the time when we took the sample, 
at least one of the lineages
in $\xi_0$, will be of type $a$.
Conversely, if $\xi_t=\emptyset$, then all lineages in $\xi_0$ are descended from
lineages that mutated
into type $A$ (having possibly first experienced merging and branching due to resampling
and selection). As a result, all individuals in $\xi_0$ are of type $A$.
We have established the following.
\begin{PRP}
Let $\xi_0=\{1,\ldots,\kappa_0\}$, then
\[\{\xi_t\neq\emptyset\} = \{|\{i: i\in\xi_0,\ X_i=a\}|\ge 1\}. \]
Hence
\[ \P(\kappa_t>0) = \E[1-(1-P_t)^{\kappa_0}]. \]
\label{prp:dual}
\end{PRP}

\subsection{Stationary distribution of the Q-process of the SDE}
\label{Q-process of diffusion}

In order to find an expression for the transition rates in our conditioned logistic
branching process we shall let $N\rightarrow\infty$ and use the distribution of the 
diffusion process~(\ref{eq:sde}).  Because we wish to study the Q-process corresponding to logistic 
branching, we shall 
need an expression for the stationary distribution of the
Q-process of~(\ref{eq:sde}). To that end, we define
\[ \beta(p) = -\mu p + sp(1-p) \]
and $\sigma(p)=\sqrt{\nsde p(1-p)}$ so that~(\ref{eq:sde}) becomes
\[ dp = \beta(p) \ dt + 
\sigma(p)\  dW, \ p_0=1. \]
We should like to find the stationary distribution of $p$ conditioned not to
hit 0. Following Karlin \& Taylor~(1981), we first calculate the speed and scale of the unconditioned
diffusion.
\nocite{karlin/taylor:1981}
Using their notation, we obtain:
\begin{eqnarray}
  s(x) &=& e^{-\frac 2 \nsde \int_0^x s-\mu/(1-\zeta) \ d\zeta}
    = e^{-2sx/\nsde} (1-x)^{-2\mu/\nsde} \label{eq:speed} \\
  S(x) &=& \int_0^x s(\zeta) \ d\zeta
    = \int_0^x e^{-2s\zeta/\nsde} (1-\zeta)^{-2\mu/\nsde} \ d\zeta \nnb \\
  m(x) &=& \frac 1 {s(x) \sigma(x)^2}
    = e^{2sx/\nsde} (1-x)^{2\mu/\nsde-1} (\nsde x)^{-1} \nnb
\end{eqnarray}
If $\mu<\nsde/2$, then $S(1)<\infty$, which means that both 0 and 1 are
hit with positive probability started from any point $p$ in $(0,1)$.
In this case, if the process $p$ is conditioned on the event $\{T_1<T_0\}$,
then according to
calculations found in Chapter 15, equations (9.1)-(9.7) of Karlin \& Taylor~(1981),
the conditioned drift is
\begin{eqnarray}
  \beta^*(x) = \beta(x) + \frac{s(x)}{S(x)} \sigma(x)^2
\label{eq:beta_star}
\end{eqnarray}
On the other hand, if $\mu\ge\nsde/2$, then $s(x)$ in~(\ref{eq:speed})
integrates to $\infty$ near 1, which makes 1 an entrance boundary.
Therefore we condition on the event $\{T_{1-\epsilon}<T_0\}$. This gives
the following conditioned drift
\[ \beta^*(x) = \left\{ \begin{array}{ll}
  \beta(x) + \frac{s(x)}{S(x)} \sigma(x)^2, & \mbox{if $x<1-\epsilon$} \\
  \beta(x), & \mbox{if $x\ge 1-\epsilon$}
  \end{array} \right. .
\]
Taking $\epsilon\to 0$ yields the same conditioned drift $\beta^*$
in~(\ref{eq:beta_star}).
Let $p^*$ denote
the process $p$ conditioned on the event $\{T_0=\infty\}$,
then the drift for $p^*$ is given by~(\ref{eq:beta_star}). We can compute
the speed and scale functions for the conditioned diffusion $p^*$.  This yields:
\begin{eqnarray*}
  s^*(x) &=& \frac{s(x)}{S(x)^2} \\
  S^*(x) &=& -\frac 1 {S(x)} \\
  \beta^*(x) &=& \beta(x) + \frac{s(x)}{S(x)} \sigma(x)^2 \\
  m^*(x) &=& \frac{1}{s^*(x)\sigma^*(x)^2} = S(x)^2 m(x) \\
  &=& \frac 1 {\nsde x} \left(\int_0^x
    e^{-2s\zeta/\nsde} (1-\zeta)^{-2\mu/\nsde} \ d\zeta\right)^2
    e^{2sx/\nsde} (1-x)^{2\mu/\nsde-1}.
\end{eqnarray*}

We now identify the stationary distribution of the conditioned process
$p^*$. We divide into two cases.
First, we consider the case $\mu<\nsde$. Since
\[ \int_0^x e^{-2s\zeta/\nsde} (1-\zeta)^{-2\mu/\nsde} \ d\zeta
  < C (1-x)^{1-2\mu/\nsde}, \]
we have
\[ m^*(x) < C(1-x)^{1-2\mu/\nsde}. \]
Hence $m^*$ is integrable near 1. Near 0, the Taylor expansion of
$\int_0^x e^{-2s\zeta/\nsde} (1-\zeta)^{-2\mu/\nsde} \ d\zeta$
is $x+\ldots$, hence $m^*$ is also integrable. This
means that the conditioned speed density $m^*$ is integrable in $[0,1]$.
Since neither 0 nor 1 is an absorbing boundary for the conditioned diffusion
$p^*$, the process is recurrent.  Therefore Theorem~4.2
of Watanabe \& Motoo~(1958) implies that $C m^*(x)$
\nocite{watanabe/motoo:1958}
is the stationary distribution of $p^*$ and, moreover, it is ergodic.

Now we deal with the case $\mu\ge\nsde$.
According to Chapter 15, equation~(5.34) of Karlin \& Taylor~(1981), the
stationary distribution $\pi^*$ of $p^*$ can be calculated from the
its speed and scale functions:
\begin{eqnarray*}
  \pi^*(x) &=& C_1 \frac{S^*(x)}{s^*(x)\sigma^*(x)^2}
    + C_2 \frac{1}{s^*(x)\sigma^*(x)^2}
  = C_1 \frac{S(x)}{s(x)\sigma(x)^2} + C_2 \frac{S(x)^2}{s(x)\sigma(x)^2} \\
  &=& m(x) (C_1 S(x) + C_2 S(x)^2).
\end{eqnarray*}
From~(\ref{eq:speed}), there exists $c>0$ such that
\begin{eqnarray}
  S(x) &\ge& c (1-x)^{1-2\mu/\nsde} \nnb \\
  m(x) &\ge& c (1-x)^{2\mu/\nsde-1}.
\label{eq:mS_lower}
\end{eqnarray}
Hence $m^*(x)=m(x) S(x)^2$ is not integrable near 1 if $\mu\ge\nsde$,
ruling out $m(x) S(x)^2$ from contributing to $\pi^*$.
Therefore the stationary distribution of $p^*$
is $\pi^*(x)=C m(x) S(x)$. We summarise the calculations above in the following
proposition.

\begin{PRP}
The diffusion process $p$ defined in~(\ref{eq:sde})
conditioned never to hit 0, denoted by $p^*$, has stationary distribution
\[ \pi^*(x) = \left\{ \begin{array}{ll}
  C m(x) S(x)^2, & \mbox{if $\mu\ge\nsde$} \\
  C m(x) S(x), & \mbox{if $\mu<\nsde$} \end{array} \right. ,\]
where $C$ is the normalising constant, and $m$ and $S$ are given
in~(\ref{eq:speed}). Moreover, $p^*$ is ergodic if $\mu<\nsde$.
\label{prp:pi_star}
\end{PRP}

\section{Logistic Branching Conditioned to Survive}
\label{transition rates}

We are now in a position to write down the transition rates for
our logistic branching process conditioned on survival and that is
our aim in \S\ref{discrete conditioning}.  In \S\ref{scaling} we record
the scaling that allows us to recover the Feller diffusion with logistic
growth and the corresponding conditioned diffusion.

\subsection{Conditioning the discrete logistic branching process}
\label{discrete conditioning}

Recall our notation, $\tau_0 = \inf\{t\ge 0: Z_t=0\}$, and let us write
$Z^T$ for the process $Z$ conditioned to survive
until time $T$, and $\P^T$ for the corresponding probability measure.
That is
\begin{eqnarray}
  \P^T(Z^T\in\cdot) = \P(Z\in\cdot \ | \ \tau_0>T)
\label{eq:PT}
\end{eqnarray}
\begin{PRP}
\label{prp:qT}
Let $\{q^T_{ij}(t)\}_{i,j\geq 1}$ be the jump rates of $Z^T$ at time $t<T$.
We have
\[ q^T_{k,k+1}(t) = q_{k,k+1} r^T_{k+1,k}(T-t), \
  q^T_{k,k-1}(t) = q_{k,k-1} r^T_{k-1,k}(T-t), \]
where
\begin{eqnarray}
  r^T_{i,j}(t) = \frac{\E[1-(1-p_t)^i]}{\E[1-(1-p_t)^j]}.
\label{eq:rT}
\end{eqnarray}
and $p$ solves the SDE~(\ref{eq:sde}).
\end{PRP}
\proof
We prove the result for $q^T_{k,k+1}$. The proof
for $q^T_{k,k-1}$ is entirely similar.

With a slight abuse of notation, let us write,
$\xi^{(k,u)}$ for the genealogical process of a sample of size $k$ taken
at time $u\leq T$ from a population evolving according to the Moran model
of \S\ref{duality}. We write $\kappa^{(k,u)}(t)$ for its size at time $t>u$.
We define
\[ \hat r_{k+1,k}(u)=\frac{\P(\kappa^{(k+1,u)}(T)>0)}{\P(\kappa^{(k,u)}(T)>0)}.\]
By Proposition~\ref{prp:dual},
\[ \hat r_{k+1,k}(u) = \frac{\E[1-(1-P_{T-u})^{k+1}]}{\E[1-(1-P_{T-u})^k]}. \]
For $a,b,a',b'>0$, we have
\[ \left|\frac a b - \frac{a'}{b'}\right|
  \le \frac a {bb'} |b'-b| + \frac 1 {b'} |a-a'|. \]
Since we can find $C_T$ such that
$\E[1-(1-P_{T-u})^k]$, $\E[1-(1-P_{T-u})^{k+1}]$,
$\E[1-(1-p_{T-u})^k]$ and $\E[1-(1-p_{T-u})^{k+1}]$ all lie in $[1/C_T,C_T]$
for $u\in [0,T]$,
we apply Proposition~\ref{prp:sde_conv} to obtain that for $u\in [0,T]$,
\begin{eqnarray}
  \left| \hat r_{k+1,k}(u) - \frac{\E[1-(1-p_{T-u})^{k+1}]}{\E[1-(1-p_{T-u})^k]}
  \right| \le C_T N^{-1/2}.
\label{eq1}
\end{eqnarray}
By an application of Bayes' rule, the jump rates of $Z^T$ at time $t$ from $k$ to $k+1$ is
\[ q^T_{k,k+1}(t) = q_{k,k+1} \frac{\P(\tau_0>T \ | \ Z_t=k+1)}{
  \P(\tau_0>T \ | \ Z_t=k)}. \]
Evidently the ratio of probabilities is uniformly bounded.  (This follows from a simple
coupling argument: take two copies of
the logistic branching process, one started from $k$ individuals and one from 
$k+1$, and wait until the first jump; with strictly positive probability, the first 
jump experienced by either of them makes the two processes equal, after which we 
use the same driving noise.)
Thanks to Lemma~\ref{lem:couple1}, the probability ratio above can be
replaced with the corresponding ratio involving $\kappa$, so that
\begin{eqnarray}
  \left|\frac{\P(\tau_0>T \ | \ q_t=k+1)}{\P(\tau_0>T \ | \ q_t=k)}
  - \frac{\P(\kappa^{(k+1,t)}(T-t)>0)}{\P(\kappa^{(k,t)}(T-t)>0)} \right|
  \le \frac{c_T}{N^{1/3}}.
\label{eq2}
\end{eqnarray}
Combining~(\ref{eq1}) and~(\ref{eq2}), we obtain
\[ q^T_{k,k+1}(t)
  = q_{k,k+1} \left( \frac{\E[1-(1-p_{T-t})^{k+1}]}{\E[1-(1-p_{T-t})^k]}
    + \Ocal(N^{-1/3})\right), \]
and the desired result follows upon taking $N\to\infty$.
\qed

We now take $T\to\infty$ in~(\ref{eq:PT}), to obtain the
Q-process, $Z^*$ (of Definition~\ref{defn:Q-process}), 
for our logistic branching process.
We denote the corresponding transition rates by $q^*_{ij}$.

~

\noindent
{\em Proof of Theorem~\ref{cor:q_star}.}

Proposition~\ref{prp:pi_star} implies that for fixed $t$,
\[ \lim_{T\to\infty} r^T_{i,j}(t) = r^*_{i,j}. \]
The formulae for $q^*_{k,k+1}$ and $q^*_{k,k-1}$ are
easy consequences of Proposition~\ref{prp:qT}.

For the bounds on $r^*_{k+1,k}$ and $r^*_{k-1,k}$, let
\[ A_k=\{\mbox{at least 1 type-a in a sample of $k$}\},\]
then
\[ \P(A_k) = \P(A_k \ | \ A_{k+1}) \P(A_{k+1}). \]
We can obtain a sample of size $k$ by removing one individual from the
sample of size $k+1$. The only way for $A_{k+1}$ to occur but
$A_k$ not to occur is if there is exactly one type-a individual in the sample
of size $k+1$ and we have removed it, therefore
\[ \P(A_k \ | \ A_{k+1}) \ge 1-\frac 1 {k+1}. \]
Hence
\[ \P(A_{k+1}) \ge \P(A_k) \ge \left(1-\frac 1 {k+1}\right) \P(A_{k+1}). \]
Since this holds regardless of the $\pi^*$, we can plug this estimate
into~(\ref{def:rstar}) to obtain
\[ 1 \le r^*_{k+1,k}
  \le \frac{\int_0^1 (1-(1-\zeta)^{k+1}) \pi^*(\zeta) \ d\zeta}{
    (1-\frac 1 {k+1}) \int_0^1 (1-(1-\zeta)^{k+1}) \pi^*(\zeta) \ d\zeta}.
\]
The desired bounds for $r^*_{k+1,k}$ and $r^*_{k-1,k}$ follow easily.
\qed

The conditioned logistic branching process $Z^*$ is a generalised birth
death process. The uniform boundedness of
$r^*_{k+1,k}$ and $r^*_{k-1,k}$ is enough to guarantee that condition (6.11.3) in
Grimmett \& Stirzaker~(1992) holds.
\nocite{grimmett/stirzaker:1992}
As a result, $Z^*$ has a unique stationary distribution and the process is
reversible under this stationary distribution.

Suppose we take a sample of size $\tilde n(0)$ from the conditioned logistic
branching process $Z^*$ at time $t$ and trace back its
ancestors. Let $\tilde n(s)$, $s\ge 0$, denote the process that counts the
number of ancestral lineages
as we trace backwards-in-time, and $\tilde z^*$ denote the time-reversed
conditioned logistic branching process.
We can write down the generator $\tilde\Gcal$ of the
process $\{\tilde z^*,\tilde n\}$ in terms of the jump rates obtained
in Theorem~\ref{cor:q_star}.

\begin{THM}
\label{thm:coalescent}
The generator $\tilde\Gcal$ of the process $\{\tilde z^*,\tilde n\}$ is given by
\begin{eqnarray*}
  \tilde\Gcal f(z,n) &=& \left(1-\frac{n(n-1)}{z(z-1)}\right) q^*_{z,z-1}
      (f(z-1,n)-f(z,n)) \\
  && + \frac{n(n-1)}{z(z-1)} q^*_{z,z-1} (f(z-1,n-1)-f(z,n)) \\
  && + q^*_{z,z+1} (f(z+1,n)-f(z,n)).
\end{eqnarray*}
\end{THM}
\proof
The rate at which $z^*$ jumps from $z$ to $z-1$ is $q^*_{z,z-1}$. Given this
occurs at time $u$, the probability it involves two lineages ancestral
to sample is ${n \choose 2}/{z \choose 2}$, hence the transition from 
$(z,n)$ to $(z-1,n-1)$.
Otherwise, the jump has no effect on the ancestral lineages to the sample,
hence the transition $(z,n)$ to $(z-1,n)$.
The transition from $z$ to $z+1$ signifies a death in the forwards-in-time
process, hence has no effect on the lineages ancestral to the sample.
\qed

We remark that the ancestral pedigree process, that is the ancestral lineages
of the entire population at present, can simply be obtained
by taking our sample $\tilde n(0)$ to the entire population.

\subsection{Rescaling the process}
\label{scaling}

In this short subsection, for completeness, we recall the rescaling of our logistic branching process 
that leads, in the limit, to a Feller diffusion with logistic growth.  
\begin{LEM}
Take
\[ s=\frac 1 2, \ \mu=\frac 1 2 - \frac b K, \ \nsde=\frac c {K^2} \]
and define \[ X = \frac 1 K Z^{1/2,c/K^2,1/2-b/K}, \]
then as $K\rightarrow\infty$, $X$ converges weakly to the solution of
\begin{equation}
\label{feller diffusion}
dX = (bX-cX^2) \ dt + \sqrt X \ dW. 
\end{equation}
\end{LEM}
The proof is standard.  
In fact one can prove a stronger result.  Using the technology of 
Barton {\em et al.}~(2004), and the calculations of the previous section,
\nocite{barton/etheridge/sturm:2004}
one can prove joint convergence under this 
scaling of the 
time reversal of our conditioned logistic branching process to the 
Q-process corresponding to~(\ref{feller diffusion}) and the genealogy of 
a sample from the population to the corresponding time-changed 
Kingman coalescent. 

\section{Populations with weak competition}
\label{sec:weak_comp}
\label{perturbation}

The difficulty with the rates established in Theorem~\ref{thm:coalescent} is that our 
expression for $r_{i,j}^*$ is difficult to evaluate, even numerically.  In this section 
we present a simple approximation of this quantity, valid for populations with only weak competition
between individuals and for which the intrinsic growth rate, $s-\mu$ is strictly positive.
This second restriction is not unreasonable from a practical perspective: its biological
interpretation is that in the absence of competition, the population would be viable. 
 
\begin{figure}[!ht]
\centering
\includegraphics[height=0.3\textheight,width=0.8\textwidth]{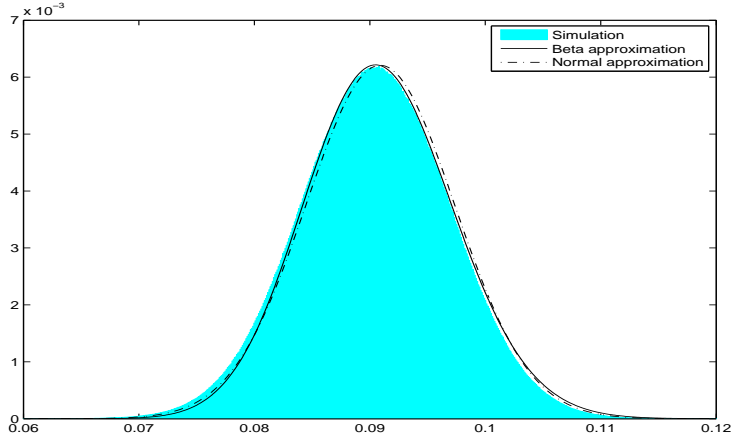}
\caption{Stationary distribution of the Q-process from simulation
versus approximation. Parameters: $s=1.1, \mu=1, \nsde=0.0001$. Bar graph
from simulation is plotted using values of solutions to the conditioned
SDE with drift parameter~(\ref{eq:beta_star}).  The beta distribution approximation
is taken from~(\ref{eq:beta}).  The normal approximation has the same mean and 
variance as those of the beta approximation.}
\label{fig:q_sde}
\end{figure}

First observe that using the approach of Norman~(1975)\nocite{norman:1975}, we see that
if competition is weak, corresponding to $\nsde$ being close to zero, then for large 
times, the solution to equation~(\ref{eq:sde})
can be approximated by that of
\begin{equation}
\label{linearised equation}
dp=s\bar{p}(\bar{p}-p)dt+\sqrt{\nsde p(1-p)}dW,
\end{equation}
where $\bar{p}=1-\mu/s$.
This is obtained by
linearising the drift in equation~(\ref{eq:sde})
around the fixed point, $\bar{p}$, of the deterministic 
equation obtained by setting $\nsde=0$.
This linearised equation is a Wright-Fisher diffusion with stationary distribution
\begin{eqnarray}
  \phi(p)=Cp^{2\alpha\bar{p}-1}(1-p)^{2\alpha(1-\bar{p})-1},
\label{eq:beta}
\end{eqnarray}
where $C$ is a normalising constant and $\alpha=2s\bar{p}/\nsde= s(s-\mu)/\nsde$.
In Figure~\ref{fig:q_sde} we examine the accuracy of this approximation.

Using $\phi$ to approximate $\pi^*$ in our expression for $r_{k+1,k}^*$ yields
$$r_{k,k+1}^*=
\frac{1-\frac{\Gamma(2\nu_1+2\nu_2)}{\Gamma(2\nu_1+2\nu_2+k+1)}\frac{\Gamma(2\nu_1+k+1)}{\Gamma(2\nu_1)}}
{1-\frac{\Gamma(2\nu_1+2\nu_2)}{\Gamma(2\nu_1+2\nu_2+k)}\frac{\Gamma(2\nu_1+k)}{\Gamma(2\nu_1)}},
$$
where 
$$\nu_1=\alpha(1-\bar{p}), \quad \nu_2=\alpha\bar{p},$$
and
$\Gamma$ is the usual gamma function.
If we consider the limit as $\nsde\downarrow 0$, using that $\Gamma(n+j)/\Gamma(n)\sim n^j$ for large $n$,
we obtain, as we should,
$$\lim_{\nsde\downarrow 0}r_{k+1,k}^*=\frac{1-\rho^{k+1}}{1-\rho^k},$$
where $\rho=\nu_1/(\nu_1+\nu_2)=\mu/s$ is the probability that a birth-death process with birth rate $s$ and death rate $\mu$
dies out in finite time.  

\section{The Yaglom limit}
\label{yaglom}

In this section, we derive (somewhat heuristically) an expression for the 
probability generating function of the Yaglom limit for our logistic branching process.
Evidently it will be dominated by the stationary distribution of the Q-process, and 
so we may take the function $g(t)$ in Definition~\ref{defn:yaglom limit} to be 
identically equal to one.

For ease of typing, we omit all superscripts in our logistic branching process.
First observe that, for $k\geq 1$,
\begin{eqnarray*}
\P(Z_{t+\epsilon}=k|Z_{t+\epsilon}\neq 0)
&=&\frac{\P(Z_{t+\epsilon}=k|Z_t\neq 0)\P(Z_t\neq 0)}{\P(Z_{t+\epsilon}\neq 0)}\\
&=&\frac{\P(Z_{t+\epsilon}=k|Z_t\neq 0)}{\P(Z_{t+\epsilon}\neq 0|Z_t\neq 0)}\\
&=&\sum_{j=1}^\infty\frac{\P(Z_{t+\epsilon}=k|Z_t=j)\P(Z_t=j|Z_t\neq 0)}{\P(Z_{t+\epsilon}\neq 0|Z_t\neq 0)}.
\end{eqnarray*}
We now estimate this quantity for small $\epsilon$.  It is convenient to write $b_k$ and
$d_k$ for the birth and death rates in the (unconditioned) logistic branching process when
the population size is $k$ and 
$\pi_t(k)=\P(Z_t=k|Z_t\neq 0)$.  Then
$$\pi_{t+\epsilon}(k)=\frac{\pi_t(k)(1-\epsilon b_k-\epsilon d_k)+\epsilon b_{k-1}\pi_t(k-1)
+\epsilon d_{k+1}\pi_t(k+1)}{1-\epsilon d_1\pi_t(1)}+{\mathcal O}(\epsilon^2).$$
Recasting this as a system of differential equations for $\pi_t(k)$ and looking for 
a fixed point, which we shall denote by $\pi(k)$, we obtain
$$b_{k-1}\pi (k-1)+d_{k+1}\pi(k+1)=(b_k+d_k-d_1\pi(1))\pi(k).$$
To obtain the probability generating function, $G(\theta)=\sum_{k=1}^\infty \pi(k)\theta^k$,
of such a fixed point, we multiply
by $\theta^k$ and sum over $k$ to obtain, on substituting for $b_k$ and $d_k$ and 
writing $a=d\pi(1)$,
\begin{eqnarray*}
0&=& \sum_{k=1}^\infty [d(k+1)+ck(k+1)]\pi(k+1)\theta^k
\\&&
-\sum_{k=1}^\infty[(b+d)k+ck(k-1)-a]
\pi(k)\theta^k+\sum_{k=1}^\infty b(k-1)\pi(k-1)\theta^k\\
&=&
\sum_{k=2}^\infty [dk+ck(k-1)]\pi(k)\theta^{k-1}
\\&&
-\sum_{k=1}^\infty[(b+d)k+ck(k-1)-a]\pi(k)\theta^k
+\sum_{k=1}^\infty bk\pi(k)\theta^{k+1}\\
&=&dG'(\theta)-d\pi(1)+c\theta G''(\theta)-(b+d)\theta G'(\theta)-c\theta^2G''(\theta)
+aG(\theta)+b\theta^2 G'(\theta).
\end{eqnarray*}
In other words, $G(\theta)$ solves
\begin{equation*}
c\theta (1-\theta)G''(\theta)+ (d-b\theta)(1-\theta))G'(\theta)+aG(\theta)=a,
\end{equation*}
where
$$ G(0)=0, \qquad G(1)=1,\quad a=dG'(0).$$

Let $D=(0,1)$ and $X$ be a process satisfying the following SDE:
\begin{eqnarray}
  dX = (d-bX)(1-X) \ dt + \sqrt{2cX(1-X)} \ dW, \quad X_0=\theta.
\label{eq:yaglomx}
\end{eqnarray}
Let $\tau_D$ be the exit time of $X$ from $D$ and
\[ f_t = e^{a (t\wedge\tau_D)} (1-G(X_{t\wedge\tau_D})) \]
By It\^o's formula, we have
\begin{eqnarray*}
  f_t &=& f_0 + \int_0^{t\wedge\tau_D} e^{a s} (a - a G(X_s)
      - (d-b X_s)(1-X_s) G'(X_s) - c X_s(1-X_s) G''(X_s) \ ds \\
  && \qquad - \int_0^{t\wedge\tau_D} a e^{a s}
      (d-b X_s) (1-X_s) G'(X_s) \ dW_s.
\end{eqnarray*}
Since the first integrand above is 0, we have
\[ \E_\theta[f_t] = 1-G(\theta). \]
Taking $t\to\infty$ yields
\[ G(\theta) = 1- \E_\theta[e^{a \tau_D} (1-G(X_{\tau_D}))]
  = 1-\E_\theta[ e^{a T_0} \ind_{T_0 < T_1}], \]
where $T_x$ is the hitting time of $x$ by $X$.
We have not used the boundary condition $a = dG'(0)$ yet. But as the Yaglom
limit is unique, which can be shown along the lines of the proof
of Theorem 8.2 of Cattiaux {\em et al}~(2009),
only one value of $a$ can satisfy the above equation.
Hence we have established the following theorem.

\begin{THM}
The probability generating function $G(\theta)$ of the Yaglom limit of
the logistic branching process defined in
Definition~\ref{defn:logistic branching} can be written as
\[ G(\theta)=1-\E_\theta[ e^{a T_0} \ind_{T_0 < T_1}], \]
where $a$ uniquely satisfies $G'(0)=a/d$
and $T_x$ is the hitting time of $x$ by the diffusion process defined
in~(\ref{eq:yaglomx}).
\end{THM}

\bibliography{structure}
\end{document}